\documentclass[12pt]{article}
\usepackage{cite}
\usepackage{mathrsfs}
\usepackage{amssymb,amsfonts,amsmath}
\usepackage{enumerate}
\usepackage{indentfirst}
\usepackage{latexsym,bm}
\usepackage[noend]{algpseudocode}
\usepackage{algorithmicx,algorithm}
\usepackage{algorithm}
\usepackage{makeidx}
\usepackage{fancybox}
\usepackage{color}
\usepackage{multicol}
\usepackage[latin1]{inputenc}
\usepackage{graphicx}
\usepackage{epstopdf}
\usepackage{pstricks,pst-node,pst-text,pst-3d}
\usepackage{tikz}
\usepackage{amsfonts,amsmath,amssymb}
\usepackage{tikz}
\usetikzlibrary{shapes,chains,scopes,positioning,backgrounds,fit,shadows,calc,arrows.meta}
\usepackage{caption}
\newtheorem{thm}{Theorem}[section]
\newtheorem{cor}[thm]{Corollary}
\newtheorem{lem}[thm]{Lemma}

\newenvironment{pf}{{\noindent \it \bf Proof:}}{{\hfill$\Box$}\\}

\def\qed{\hfill \nopagebreak\rule{5pt}{8pt}}

\begin{document}

\title{\bf Perfect Out-forests and Steiner Cycle Packing in Digraphs}
\author{
Yuefang Sun\\
School of Mathematics and Statistics, Ningbo University,\\
Ningbo 315211, P. R. China\\ 
Email address: sunyuefang@nbu.edu.cn}
\maketitle

\begin{abstract}
In this paper, we study the complexity of two types of digraph packing problems: perfect out-forests problem and Steiner cycle packing problem.

For the perfect out-forests problem, we prove that it is NP-hard to decide whether a given strong digraph contains a 1-perfect out-forest. However, when restricted to a semicomplete digraph $D$, the problem of deciding whether $D$ contains an $i$-perfect out-forest becomes polynomial-time solvable, where $i\in \{0,1\}$. We also prove that it is NP-hard to find a 0-perfect out-forest of maximum size in a connected acyclic digraph, and it is NP-hard to find a 1-perfect out-forest of maximum size in a connected digraph.

For the Steiner cycle packing problem, when both $k\geq 2, \ell\geq 1$ are fixed integers, we show that the problem of deciding whether there are at least $\ell$ internally disjoint directed $S$-Steiner cycles in an Eulerian digraph $D$ is NP-complete, where $S\subseteq V(D)$ and $|S|=k$. However, when we consider the class of symmetric digraphs, the problem becomes polynomial-time solvable. We also show that the problem of deciding whether there are at least $\ell$ arc-disjoint directed $S$-Steiner cycles in a given digraph $D$ is NP-complete, where $S\subseteq V(D)$ and $|S|=k$.

\vspace{0.3cm}

\vspace{0.3cm}{\bf AMS subject
classification (2020)}: 05C05, 05C20, 05C38, 05C45, 05C70, 05C85, 68Q25, 68R10.

\end{abstract}

\section{Introduction}

\subsection{Motivations}

We refer the readers to \cite{Bang-Jensen-Gutin} for graph-theoretical notation and terminology not given here. Note that all digraphs considered in this paper have no parallel arcs or loops. We use $[n]$ to denote the set of all natural numbers from 1 to $n$. For a digraph $D$, we use $UG(D)$ to denote the underlying (undirected) graph of $D$. A digraph $D$ is called {\em connected} if $UG(D)$ is connected. Moreover, $D$ is {\em strong connected} or, simply, {\em strong}, if for any pair of vertices $x, y\in V(D)$, there is a path from $x$ to $y$ in $D$, and vice versa.

A digraph is {\em acyclic} if it has no directed cycle. An {\em out-tree} (resp. {\em in-tree}) {\em rooted at a vertex $r$} is an orientation of a tree such that the in-degree (resp. out-degree) of every vertex but $r$ equals one. An {\em out-branching} $B^+_r$ (resp. {\em in-branching} $B^-_r$) in a digraph $D$ is a spanning subdigraph of $D$ which is an out-tree (resp. in-tree). 

For $i \in \{0, 1\}$, a spanning forest $F$ of a graph $G$ is an {\em $i$-perfect forest} if each tree of $F$ is an induced subgraph of $G$, and exactly $i$ vertices in $G$ have even degree (including zero). A 0-perfect forest is also called a perfect forest. Clearly, the concept of perfect forest is a natural generalization of perfect matching of graphs. The problem of $i$-perfect forest in undirected graphs has been studied by some researchers \cite{Caro-Lauri-Zarb, GutinJGT2016, Gutin-YeoMFCS2021, ScottGC2001, Sharan-Wigderson1996}. 

Gutin and Yeo \cite{Gutin-YeoJGT2017} introduced the concept of perfect out-forest in digraphs.
An {\em out-forest} is a collection of vertex disjoint out-trees. A spanning out-forest $F$ of a digraph $D$ is a {\em perfect out-forest} if each out-tree of $F$ is an induced subgraph of $D$, and the degree of each vertex in $UG(F)$ is odd. Clearly, the concept of perfect out-forest is a natural generalization of perfect matching of digraphs. Similarly, we define the concept of $i$-perfect out-forest as follows. 
For $i \in \{0, 1\}$, a spanning out-forest $F$ of a digraph $D$ is an {\em $i$-perfect out-forest} if each out-tree of $F$ is an induced subgraph of $D$, and exactly $i$ vertices in $UG(F)$ have even degree (including zero). Clearly, a 0-perfect out-forest is exactly a perfect out-forest.

For a graph $G=(V(G),E(G))$ and a set $S\subseteq V(G)$ of at
least two vertices, an {\em $S$-Steiner tree} or, simply, an {\em
$S$-tree} is a tree $T$ of $G$ with $S\subseteq V(T)$. Two $S$-trees are said to be {\em edge-disjoint} if they have no common edge. Two edge-disjoint $S$-trees $T_1$ and $T_2$ are said to be {\em internally disjoint} if $V(T_1)\cap V(T_2)=S$. Let $D=(V(D),A(D))$ be a digraph with order $n$. Let $S\subseteq V(D)$ with
$r\in S$ and $2\leq |S|\leq n$. A {\em directed $(S, r)$-Steiner tree} or, simply, an {\em $(S, r)$-trees}, is an out-tree $T$ rooted at $r$ with
$S\subseteq V(T)$. Two $(S, r)$-trees
are said to be {\em arc-disjoint} if they have no common arc. Two arc-disjoint $(S, r)$-trees $T_1$ and $T_2$ are
said to be {\em internally disjoint} if $V(T_1)\cap V(T_2)=S$. A strong subgraph $H$ of $D$ is called an {\em $S$-strong
subgraph} if $S\subseteq V(H)$. Two $S$-strong subgraphs are said to be {\em arc-disjoint} if they have no common arc. Furthermore, two arc-disjoint $S$-strong subgraphs $D_1$ and $D_2$ are said {\em internally disjoint} if $V(D_1)\cap V(D_2)=S$.

The basic problem of {\sc Steiner Tree
Packing} (also called {\sc edge-disjoint Steiner tree packing problem}) is defined as follows: the input consists of an undirected
graph $G$ and a subset of vertices $S\subseteq V(D)$, the goal is to
find a largest collection of edge-disjoint $S$-Steiner trees.  The Steiner tree packing problem has applications in VLSI circuit design \cite{Grotschel-Martin-Weismantel, Sherwani}. 
In this application, a Steiner tree is needed to share an electronic signal by a set of terminal nodes. Besides the classical version, people also study some other
variations, such as internally disjoint Steiner tree packing problem \cite{Cheriyan-Salavatipour,
Li-Mao5}.

There are two ways to extend the problem of Steiner tree packing to directed graphs: directed Steiner tree packing problem \cite{Cheriyan-Salavatipour, Sun-Yeo} and strong subgraph packing problem \cite{Sun-Gutin2, Sun-Gutin-Yeo-Zhang}. In the problem of {\sc Arc-disjoint directed Steiner tree packing} (resp. {\sc Internally disjoint directed Steiner tree packing}), one aims to find a largest collection of arc-disjoint (resp. internally disjoint) directed $(S,r)$-Steiner trees in $D$, while in the problem of {\sc Arc-disjoint strong subgraph packing} (resp. {\sc Internally disjoint strong subgraph packing}), the goal is to find a largest collection of arc-disjoint (resp. internally disjoint) $S$-strong subgraphs, where $S\subseteq V(D)$ and $r\in S$. The directed Steiner tree packing problem and strong subgraph packing problem both belong to Steiner type packing problems in digraphs.

In this paper, we introduce a new type of Steiner type packing problem in digraphs. Let $D=(V(D),A(D))$ be a digraph of order $n$, $S\subseteq V(D)$ a $k$-subset of $V(D)$ and $2\le k\leq n$. A directed cycle $C$ of $D$ is called a {\em directed $S$-Steiner cycle} or, simply, an {\em $S$-cycle} if $S\subseteq V(C)$. Two $S$-cycles are said to be {\em arc-disjoint} if they have no common arc. Furthermore, two arc-disjoint $S$-cycles are {\em internally disjoint} if the set of common vertices of them is exactly $S$. 
We define the following problems of packing Steiner cycles in digraphs as follows:

{\sc Arc-disjoint directed Steiner cycle packing:}
The input consists of a digraph $D$ and a subset of vertices $S\subseteq V(D)$, the goal is to find a largest collection of arc-disjoint $S$-cycles.

{\sc Internally disjoint directed Steiner cycle packing:}
The input consists of a digraph $D$ and a subset of vertices $S\subseteq V(D)$, the goal is to find a largest collection of internally disjoint $S$-cycles.

Let $\kappa^c_S(D)$~(resp. $\lambda^c_S(D)$), be the maximum number of internally disjoint (resp. arc-disjoint) $S$-cycles in $D$. Although an $S$-cycle is also an $S$-strong subgraph, the problem of Steiner cycle packing is quite distinct with strong subgraph packing. For example, observe that it can be decided in polynomial-time whether $\kappa_S(D)\geq 1$, where $\kappa_S(D)$ denotes the maximum number of internally disjoint $S$-strong subgraphs in $D$. However, it is NP-complete to decide whether $\kappa^c_S(D)\geq 1$ even restricted to an Eulerian digraph $D$ by Theorem~\ref{thm1a}.

The problem of packing directed Steiner cycle 
is also related to other problems in graph theory. For example, when $|S|=n$, an $S$-cycle is a Hamiltonian cycle. Therefore, in this case, the problem is equivalent to finding maximum number of arc-disjoint Hamiltonian cycles. A digraph $D$ is {\em $k$-cyclic} if it has a cycle containing the vertices  $x_1, x_2, \dots, x_k$ for every choice of $k$ vertices. Note that the notion of $k$-cyclic attracts the attention of some researchers, such as \cite{Kuhn2008}. By definition, a digraph is $k$-cyclic if and only if $\kappa^c_S(D)\geq 1$ (resp. $\lambda^c_S(D)\geq 1$) for every $k$-subset $S$ of $V(D)$.


\subsection{Additional terminology and notation}

A digraph $D$ is {\em semicomplete} if for every distinct $x,y\in V(D)$ at least one of the arcs $xy,yx$ is in $D$.  A {\em tournament} is a semicomplete digraph without 2-cycles. 
A digraph $D$ is {\em symmetric} if every arc in $D$ belongs to a $2$-cycle.  
In other words, a symmetric digraph $D$ can be obtained from its underlying undirected graph $G$ by replacing each edge of $G$ with the corresponding arcs of both directions, that is, $D=\overleftrightarrow{G}.$ 
A digraph $D$ is {\em Eulerian} if $D$ is connected and $d^+(x)=d^-(x)$ for every vertex $x\in V(D)$.

\subsection{Our results}

In Section 2, by using the NP-hardness result (obtained by Gutin \& Yeo in \cite{Gutin-YeoJGT2017})  for the problem of perfect out-forest in strong digraphs,  we prove the NP-hardness of deciding whether a given strong digraph contains a 1-perfect out-forest (Theorem~\ref{thmA}(a)). However, when restricted to semicomplete digraphs, the problem of deciding whether $D$ contains an $i$-perfect out-forest becomes polynomial-time solvable, where $i\in \{0,1\}$ (Theorem~\ref{thmA}(b)). The problem of finding a perfect out-forest of size at least $n-1$ is polynomial-time solvable because $D$ has a perfect forest of size at least $n-1$ if
and only if $D$ is an out-branching in which every vertex is of odd degree in $UG(D)$. However, as shown in Theorem~\ref{thmB}, when ``$n-1$'' is replaced by ``$n-k$'' ($k\geq 2$), the problem above becomes NP-hard, even restricted to connected acyclic digraphs. Similarly, we also prove that it is NP-hard to find a 1-perfect out-forest of size at least $n-k$ for a fixed integer $k\geq 2$ (Theorem~\ref{thmC}). These two theorems imply that it is NP-hard to find a perfect out-forest of maximum size in a connected acyclic digraph, and it is NP-hard to find a 1-perfect out-forest of maximum size in a connected digraph (Theorem~\ref{thmD}).

In Section 3, we study the problem of Steiner cycle packing. When both $k\geq 2, \ell\geq 1$ are fixed integers, we show that the problem of deciding whether $\kappa^c_S(D)\geq \ell$ with $|S|=k$ is NP-complete, even restricted to Eulerian digraphs (Theorem~\ref{thm1a}). However, when we consider the class of symmetric digraphs, the problem becomes polynomial-time solvable (Theorem~\ref{thm1c}). We also show that the problem of deciding whether $\lambda^c_S(D)\geq \ell$ with $|S|=k$ is NP-complete, when both $k\geq 2, \ell\geq 1$ are fixed integers (Theorem~\ref{thm1d}).

\section{Perfect out-forests}

Gutin and Yeo proved the complexity for the problem of deciding whether a given strong digraph $D$ contains a perfect out-forest.

\begin{thm}\label{thm01}\cite{Gutin-YeoJGT2017}
The problem of deciding whether a given strong digraph $D$ contains a perfect out-forest is NP-hard.
\end{thm}

By Theorem~\ref{thm01}, we can prove a similar result for the problem of 1-perfect out-forests. These two results mean that the problem of $i$-perfect out-forests ($i\in \{0,1\}$) in a general digraph is difficult. However, when restricted to some digraph classes, such as semicomplete digraphs, the problem becomes polynomial-time solvable. Furthermore, for a tournament $T$ with order $n$, one can prove that $T$ contains a perfect out-forest if and only if $n$ is even, and $T$  contains a 1-perfect out-forest if and only if $n$ is odd.

\begin{thm}\label{thmA}
The following assertions hold:
\begin{description}
\item[(a) ] The problem of deciding whether a given strong digraph $D$ contains a 1-perfect out-forest is NP-hard.
\item[(b) ] The problem of deciding whether a given semicomplete digraph $D$ contains an $i$-perfect out-forest is polynomial-time solvable, where $i\in \{0,1\}$.
\end{description}
\end{thm}
\begin{pf}

\textbf{Part (a):} Let $D$ be any strong digraph, we construct a strong digraph $D'$ from $D$ by adding a new vertex $v$ and two arcs $uv, vu$ for some vertex $u\in V(D)$. If $D$ contains a perfect out-forest $F$, then $F'=F\cup \{v\}$ is clearly a 1-perfect out-forest of $D'$ such that $v$ is the unique vertex whose degree in the underlying graph of $F'$ is even. If $D'$ contains a 1-perfect out-forest $F'$, then $\{v\}\in F'$. Indeed, if $\{v\}\not\in F'$, then there is an out-tree $T\in F'$ containing $u, v$. As $T$ is an induced subgraph of $D$, we have $uv, vu\in A(T)$, a contradiction. Hence $\{v\}\in F'$ and $v$ is the unique vertex whose degree in the underlying graph of $F'$ is even. Now observe that $F=F'\setminus \{v\}$ is a perfect out-forest of $D$. Therefore, $D$ contains a perfect out-forest if and only if $D'$ contains a 1-perfect out-forest, by Theorem~\ref{thm01}, the result holds.

\textbf{Part (b):} Let $D$ be a semicomplete digraph. Observe that $F$ is a perfect out-forest of $D$ if and only if $F$ is a perfect matching of $D$ such that if $xy\in A(F)$, then $yx\not \in A(D)$. Now we just remove all arcs that belong to 2-cycles in $D$. That is, if $xy$ and $yx$ are both arcs in $D$, then neither can be used in the matching $F$, so delete both from the digraph. Then the problem is equivalent to deciding if the resulting digraph (with no 2-cycles), denoted by $D'$, has a perfect matching. Clearly, this problem is polynomial-time solvable, as $D'$ has a perfect matching if and only if $UG(D')$ has a perfect matching, and the latter can be decided in polynomial-time. Hence, the problem of deciding whether a given semicomplete digraph $D$ contains a perfect out-forest is polynomial-time solvable.

Observe that $D$ has a 1-perfect out-forest if and only if for some vertex $u\in V(D)$, $D-u$ has a perfect matching $F'$ such that if $xy\in A(F')$, then $yx\not \in A(D-u)$. By the argument above, we can in polynomial-time decide if $D-u$ has such a perfect matching. Hence, the problem of deciding whether a given semicomplete digraph $D$ contains a 1-perfect out-forest is polynomial-time solvable.
\end{pf}







Recall that the problem of finding a perfect out-forest of size at least $n-1$ is polynomial-time solvable because $D$ has a perfect forest of size at least $n-1$ if
and only if $D$ is an out-branching in which every vertex is of odd degree in $UG(D)$. However, as shown in Theorem~\ref{thmB}, when ``$n-1$'' is replaced by ``$n-k$'' ($k\geq 2$), the problem above becomes NP-hard.

The {\em not-all-equal 3-SAT} problem, abbreviated to NAE-3-SAT, is the problem of determining whether an instance of 3-SAT has a truth assignment to its variables such that each clause contains both a true and a false literal. If such case holds for an instance $I$, then we say that $I$ is {\em NAE-satisfied}.

\begin{thm}\label{thmB}
Let $k\geq 2$ be an integer. The problem of deciding whether a given connected acyclic digraph $D$ of order $n$ contains a perfect out-forest with at least $n-k$ arcs is NP-hard.
\end{thm}
\begin{pf} We prove the result by induction on $k$. We first consider the base step that $k=2$, that is, we show that the problem of deciding whether a given connected acyclic digraph $D$ of order $n$ contains a perfect out-forest with at least $n-2$ arcs is NP-hard, by reducing from the NAE-3-SAT problem which is NP-hard \cite{Schaefer}. 

\tikzset{arrow1/.style = {draw = black, thick, -{Latex[length = 3.5mm, width = 1.1mm]},}
	}
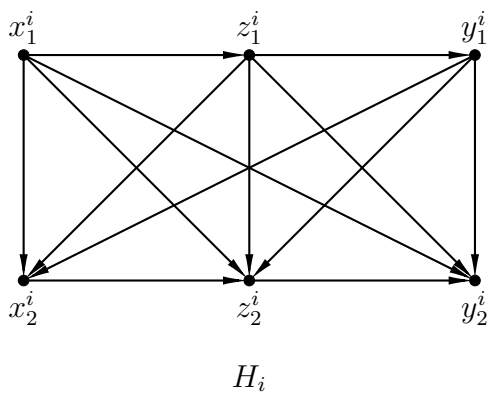
\begin{figure}[htb]
		\centering
		\begin{tikzpicture}
			\filldraw[black]    (0, 0)  circle (2pt)  node [anchor=south] {$x^{i}_{1}$};
			\filldraw[black]    (3, 0)  circle (2pt)  node [anchor=south] {$z^{i}_{1}$};
			\filldraw[black]    (6, 0)  circle (2pt)  node [anchor=south] {$y^{i}_{1}$};
			\filldraw[black]    (0, -3)  circle (2pt)  node [anchor=north] {$x^{i}_{2}$};
			\filldraw[black]    (3, -3)  circle (2pt)  node [anchor=north] {$z^{i}_{2}$};
			\filldraw[black]    (6, -3)  circle (2pt)  node [anchor=north] {$y^{i}_{2}$};
			
			\draw[arrow1] []      (0, 0) -- (3, 0);
			\draw[arrow1] []      (3, 0) -- (6, 0);
			\draw[arrow1] []      (0, 0) -- (0, -3);
			\draw[arrow1] []      (0, 0) -- (3, -3);
			\draw[arrow1] []      (0, 0) -- (6, -3);
			\draw[arrow1] []      (3, 0) -- (0, -3);
			\draw[arrow1] []      (3, 0) -- (3, -3);
			\draw[arrow1] []      (3, 0) -- (6, -3);
			\draw[arrow1] []      (6, 0) -- (0, -3);
			\draw[arrow1] []      (6, 0) -- (3, -3);
			\draw[arrow1] []      (6, 0) -- (6, -3);
			\draw[arrow1] []      (0, -3) -- (3, -3);
			\draw[arrow1] []      (3, -3) -- (6, -3);
			
			\node at(3, -4.3) {$H_{i}$};
		\end{tikzpicture}
		\caption{The digraph $H_i$.}
		\label{fig2}
	\end{figure}		

Let $I$ be an instance of NAE-3-SAT with clauses $C_1, C_2, \dots, C_m$ and variables $v_1, v_2, \dots, v_n$. Our construction is based on a construction from \cite{Gutin-YeoMFCS2021}. For $i\in [n]$, let $H_i$
be a digraph (see Figure~\ref{fig2}) with vertex set $$\{x^i_1, z^i_1, y^i_1, x^i_2, z^i_2, y^i_2\}$$ and arc set $$\{x^i_1z^i_1,z^i_1y^i_1,x^i_2z^i_2,z^i_2y^i_2,x^i_1x^i_2,x^i_1z^i_2,x^i_1y^i_2,z^i_1x^i_2,z^i_1z^i_2,z^i_1y^i_2,y^i_1x^i_2,y^i_1z^i_2,y^i_1y^i_2\}.$$
For $i\in [n-1]$, we add all arcs from $\{y^i_1, y^i_2\}$ to $\{x^{i+1}_1, x^{i+1}_2\}$:
$y^i_1x^{i+1}_1$, $y^i_2x^{i+1}_1$, $y^i_1x^{i+1}_2$, $y^i_2x^{i+1}_2$. For each vertex $u\in V(H_i)\setminus \{x^1_1,x^1_2,y^n_1,y^n_2\}$ ($i\in [n]$), add a new vertex $u'$ and a new arc $uu'$. Now the current digraph is denoted by $Q$. 

Let $V(D)=V(Q)\cup\{c_j,c'_j\mid j\in[m]\}$. For each $j\in[m]$, we add the arcs $y^i_2c_j,y^i_2c'_j$ (resp. $y^i_1c_j,y^i_1c'_j$) if and only if $v_i$ (resp. $\overline{v_i}$) is a literal in the clause $C_j$. This completes the construction of $D$. Observe that $D$ is a connected acyclic digraph. See Figure~\ref{fig3} for an example of $D$ when $I=(v_1,\overline{v_2},v_3)$.

\tikzset{arrow1/.style = {draw = green, thick, -{Latex[length = 3.5mm, width = 1.3mm]},}
	}
    \tikzset{arrow2/.style = {draw = red, thick, -{Latex[length = 3.5mm, width = 1.3mm]},}
    }
    \tikzset{arrow3/.style = {draw = blue, thick, -{Latex[length = 3.5mm, width = 1.3mm]},}
    }
	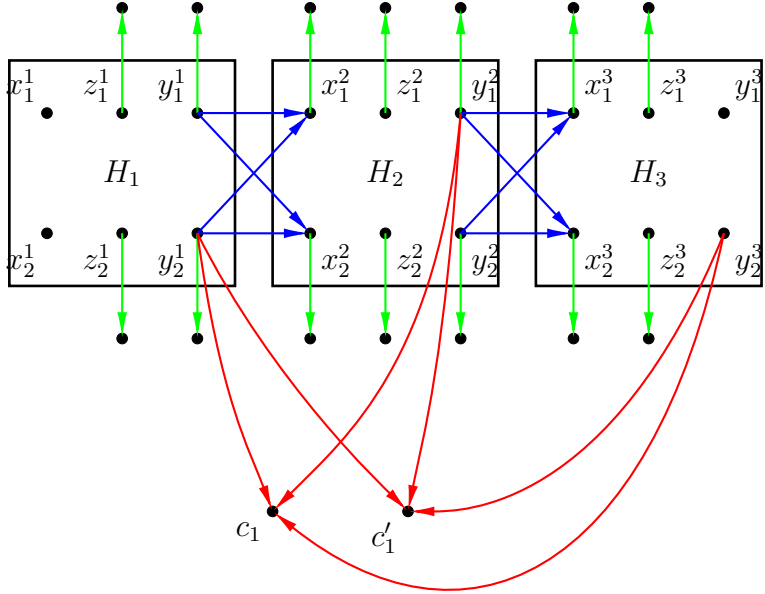
\begin{figure}[htb]
		\centering
		\begin{tikzpicture}
			\draw[line width=1pt] (0, 0) rectangle (3, -3);
			\draw[line width=1pt] (3.5, 0) rectangle (6.5, -3);
			\draw[line width=1pt] (7, 0) rectangle (10, -3);
			
			\filldraw[black]    (0.5, -0.7)  circle (2pt)  node [anchor=south east] {$x^{1}_{1}$};
			\filldraw[black]    (1.5, -0.7)  circle (2pt)  node [anchor=south east] {$z^{1}_{1}$};
			\filldraw[black]    (2.5, -0.7)  circle (2pt)  node [anchor=south east] {$y^{1}_{1}$};
			\filldraw[black]    (1.5, 0.7)  circle (2pt);
			\filldraw[black]    (2.5, 0.7)  circle (2pt);
			\filldraw[black]    (0.5, -2.3)  circle (2pt)  node [anchor=north east] {$x^{1}_{2}$};
			\filldraw[black]    (1.5, -2.3)  circle (2pt)  node [anchor=north east] {$z^{1}_{2}$};
			\filldraw[black]    (2.5, -2.3)  circle (2pt)  node [anchor=north east] {$y^{1}_{2}$};
			\filldraw[black]    (1.5, -3.7)  circle (2pt);
			\filldraw[black]    (2.5, -3.7)  circle (2pt);
	
			\filldraw[black]    (4, -0.7)  circle (2pt)  node [anchor=south west] {$x^{2}_{1}$};
			\filldraw[black]    (5, -0.7)  circle (2pt)  node [anchor=south west] {$z^{2}_{1}$};
			\filldraw[black]    (6, -0.7)  circle (2pt)  node [anchor=south west] {$y^{2}_{1}$};
			\filldraw[black]    (4, 0.7)  circle (2pt);
			\filldraw[black]    (5, 0.7)  circle (2pt);
			\filldraw[black]    (6, 0.7)  circle (2pt);
			\filldraw[black]    (4, -2.3)  circle (2pt)  node [anchor=north west] {$x^{2}_{2}$};
			\filldraw[black]    (5, -2.3)  circle (2pt)  node [anchor=north west] {$z^{2}_{2}$};
			\filldraw[black]    (6, -2.3)  circle (2pt)  node [anchor=north west] {$y^{2}_{2}$};
			\filldraw[black]    (4, -3.7)  circle (2pt);
			\filldraw[black]    (5, -3.7)  circle (2pt);
			\filldraw[black]    (6, -3.7)  circle (2pt);
			
			\filldraw[black]    (7.5, -0.7)  circle (2pt)  node [anchor=south west] {$x^{3}_{1}$};
			\filldraw[black]    (8.5, -0.7)  circle (2pt)  node [anchor=south west] {$z^{3}_{1}$};
			\filldraw[black]    (9.5, -0.7)  circle (2pt)  node [anchor=south west] {$y^{3}_{1}$};
			\filldraw[black]    (7.5, 0.7)  circle (2pt);
			\filldraw[black]    (8.5, 0.7)  circle (2pt);
			\filldraw[black]    (7.5, -2.3)  circle (2pt)  node [anchor=north west] {$x^{3}_{2}$};
			\filldraw[black]    (8.5, -2.3)  circle (2pt)  node [anchor=north west] {$z^{3}_{2}$};
			\filldraw[black]    (9.5, -2.3)  circle (2pt)  node [anchor=north west] {$y^{3}_{2}$};
			\filldraw[black]    (7.5, -3.7)  circle (2pt);
			\filldraw[black]    (8.5, -3.7)  circle (2pt);
			
			\filldraw[black]    (3.5, -6)  circle (2pt)  node [anchor=north east] {$c_{1}$};
			\filldraw[black]    (5.3, -6)  circle (2pt)  node [anchor=north east] {$c'_{1}$};
			
			\draw[arrow1] []      (1.5, -0.7) -- (1.5, 0.7);
			\draw[arrow1] []      (2.5, -0.7) -- (2.5, 0.7);
			\draw[arrow1] []      (1.5, -2.3) -- (1.5, -3.7);
			\draw[arrow1] []      (2.5, -2.3) -- (2.5, -3.7);
			\draw[arrow1] []      (4, -0.7) -- (4, 0.7);
			\draw[arrow1] []      (5, -0.7) -- (5, 0.7);
			\draw[arrow1] []      (6, -0.7) -- (6, 0.7);
			\draw[arrow1] []      (4, -2.3) -- (4, -3.7);
			\draw[arrow1] []      (5, -2.3) -- (5, -3.7);
			\draw[arrow1] []      (6, -2.3) -- (6, -3.7);
			\draw[arrow3] []      (2.5, -0.7) -- (4, -0.7);
			\draw[arrow3] []      (2.5, -0.7) -- (4, -2.3);
			\draw[arrow3] []      (2.5, -2.3) -- (4, -0.7);
			\draw[arrow3] []      (2.5, -2.3) -- (4, -2.3);
			
			\draw[arrow1] []      (7.5, -0.7) -- (7.5, 0.7);
			\draw[arrow1] []      (8.5, -0.7) -- (8.5, 0.7);
			\draw[arrow1] []      (7.5, -2.3) -- (7.5, -3.7);
			\draw[arrow1] []      (8.5, -2.3) -- (8.5, -3.7);
			\draw[arrow3] []      (6, -0.7) -- (7.5, -0.7);
			\draw[arrow3] []      (6, -0.7) -- (7.5, -2.3);
			\draw[arrow3] []      (6, -2.3) -- (7.5, -0.7);
			\draw[arrow3] []      (6, -2.3) -- (7.5, -2.3);
			
			\draw[arrow2] []      (2.5, -2.3) .. controls (2.8, -4.5) and (3.2, -5.2) .. (3.5, -6);
			\draw[arrow2] []      (2.5, -2.3) .. controls (3.5, -4.3) and (5, -5.7) .. (5.3, -6);
			\draw[arrow2] []      (6, -0.7) .. controls (5.7, -4) and (4.5, -5) .. (3.5, -6);
			\draw[arrow2] []      (6, -0.7) .. controls (5.8, -4) and (5.5, -5) .. (5.3, -6);
			\draw[arrow2] []      (9.5, -2.3) .. controls (8.5, -7) and (6, -8) .. (3.5, -6);
			\draw[arrow2] []      (9.5, -2.3) .. controls (8, -6) and (6, -6) .. (5.3, -6);
			
			\node at(1.5, -1.5) {$H_{1}$};
			\node at(5, -1.5) {$H_{2}$};
			\node at(8.5, -1.5) {$H_{3}$};
		\end{tikzpicture}
		\caption{The digraph $D$ if $I=(v_1,\overline{v_2},v_3)$.}
		\label{fig3}
	\end{figure}

To complete the proof of the base step, it suffices to show that $D$ contains a perfect out-forest of size at least $n-2$ if and only if $I$ is NAE-satisfied. 

We first assume that $D$ contains a perfect out-forest $F$ with at least $|V(D)|-2$ arcs. Observe that now $F$ consists of two vertex-disjoint out-trees $T_1$ and $T_2$ such that each $T_i$ is an induced subdigraph in $D$ of order at least two, the degree of each vertex in $UG(T_i)$ is odd, and $V(T_1)\cup V(T_2)=V(D)$. We will prove the following three claims where Claim 3 completes the proof this direction.

\textbf{Claim 1:} For each $i\in [n]$, exactly one of the following assertions hold:

A: $\{x^i_1, z^i_1, y^i_1\} \subseteq V(T_1)$ and $\{x^i_2, z^i_2, y^i_2 \} \subseteq V(T_2)$.

B: $\{x^i_1, z^i_1, y^i_1\} \subseteq V(T_2)$ and $\{x^i_2, z^i_2, y^i_2 \} \subseteq V(T_1)$.

\textbf{Proof of Claim 1:} Observe that $|V(T_j)\cap V(H_i)|=3$ for $j\in [2]$ and $i\in [n]$. Indeed, if $|V(T_j)\cap V(H_i)|\geq 4$ for some $j\in [2]$ and $i\in [n]$, then the underlying graph of $D[V(T_j)]$ contains a 4-cycle, a contradiction. 
Since there is no 3-cycle in the underlying graph of $D[V(T_j)]$, $x^i_1$ and $y^i_1$ must belong to the same tree, say $T_j$, and $x^i_2$ and $y^i_2$ must belong to the other tree $T_{3-j}$. By the construction of $Q$, exactly one of A and B holds. This completes the proof of Claim 1. \qed

\textbf{Claim 2:} $|N_D(c_j)\cap V(T_i)|=1$ whenever $c_j\in V(T_i)$.

\textbf{Proof of Claim 2:} Recall that both $T_1$ and $T_2$ have orders at least two, any vertex having degree one in $UG(D)$ must belong to the same tree as its neighbour. By Claim 1, $D[V(Q)\cap V(T_i)]$ is a directed path of length $3n$ with an pendent out-arc from each inner vertex of the path.

As both $T_1$ and $T_2$ have orders at least two, we have $|N_D(c_j)\cap V(T_i)|\geq 1$. Suppose that $|N_D(c_j)\cap V(T_i)|\geq 2$. As $T_i$ is an induced out-tree in $Q$, $c_j$ must have at least two neighbours, say $x$ and $y$, in $T_i$. By the argument above, there is an $x-y$ path in the underlying graph of $D[V(Q)\cap V(T_i)]$, then together with the vertex $c_j$, there is a cycle in $UG(T_i)$, a contradiction. This completes the proof of Claim 2. \qed



\textbf{Claim 3:} The instance $I$ is NAE-satisfiable.


\textbf{Proof of Claim 3:} We give a truth assignment as follows: if $y_1^i\in V(T_1)$, then let $v_i$ be true, and if $y_2^i\in V(T_1)$ (which means that $y_1^i\not\in V(T_1)$), then let $v_i$ be false.

Without loss of generality, we assume that $c_j\in V(T_1)$ (the argument for the case that $c_j\in V(T_2)$ is similar and hence we omit the details), and the two variables $v_{i_1}$, $v_{i_2}$ appear in $C_j$. 
We first consider the case that $y_1^{i_1}\in V(T_1)$. Now $c_jy_1^{i_1}\in E(T_1)$, so $\overline{v_{i_1}}$ is a literal of $C_j$ (by the construction of $D$). As now $v_{i_1}$ is true, the literal $\overline{v_{i_1}}$ is false. By Claim~2, we have $|N_{D}(c_j)\cap V(T_1)|=1$, so $c_jy_1^{i_2}\in A(T_2)$ or $c_jy_2^{i_2}\in A(T_2)$. For the former case, $\overline{v_{i_2}}$ is a literal of $C_j$ and $y_2^{i_2}\in V(T_1)$ (and so $v_{i_2}$ is false), then $\overline{v_{i_2}}$ is a true literal of $C_j$. For the latter case, $v_{i_2}$ is a literal of $C_j$ and $y_1^{i_2}\in V(T_1)$ (and so $v_{i_2}$ is true), then $v_{i_2}$ is still a true literal of $C_j$. 
Hence, we have that $C_j$ has at least one true literal and at least one false literal in this case. For the case that $y_2^{i_1}\in V(T_1)$, the argument is similar, so we omit the details. This completes the proof of Claim 3 and the first direction. \qed



Conversely, assume that $I$ is NAE-satisfied and consider a truth assignment $\tau$ NAE-satisfying $I$. We now construct two vertex-disjoint induced out-trees, $T_1$ and $T_2$, in $D$ such that $V(T_1)\cup V(T_2)=V(D)$, and all degrees in $UG(T_i)$ are odd for $i\in [2]$. If $v_i$ is true (resp. false) in $\tau$, then add vertices $x^i_1, z^i_1, y^i_1$ (resp. $x^i_2, z^i_2, y^i_2$) to $T_1$ and add vertices $x^i_2, z^i_2, y^i_2$ (resp. $x^i_1, z^i_1, y^i_1$) to $T_2$. We then update $T_i$ by adding $u'$ (if exists) for each vertex $u$ of the current $T_i$ ($i\in [2]$). Note that now each $V(T_i)$ induces an out-tree in $D$, where each vertex has odd degree in $UG(T_i)$. We still need the following claim.

\textbf{Claim 4:} Each of $c_j$ and $c'_j$ has exactly one arc coming from one tree and two arcs coming from the other tree, where $j\in [m]$.

\textbf{Proof of Claim 4:}
Suppose that all three arcs incident to $c_j$ (resp. $c'_j$) come from $T_1$, and the variables $v_{i_1}, v_{i_2}, v_{i_3}$ appear in $C_j$. Furthermore, without loss of generality, assume that $y^{i_1}_1, y^{i_2}_1, y^{i_2}_1 \in V(T_1)$. This means that $C_j$ consists of the following three literals: $\overline{v_{i_1}}, \overline{v_{i_2}}, \overline{v_{i_3}}$, and all of $v_{i_1}, v_{i_2}, v_{i_3}$ are true, which produces a contradiction as now $C_j$ contains no true literal. This completes the proof of Claim 4. \qed

In fact, we can further prove that the tree which has exactly one arc into $c_j$ is the same as the tree which has exactly one arc into $c'_j$ by the construction of $D$. Now add each of $c_j$ and $c'_j$ to the tree from which there is exactly one arc. It can be checked that $T_1$  and $T_2$ form a perfect out-forest in $D$ with $|V(D)|-2$ arcs, where $T_i$ is an induced out-tree rooted at $x_i^1$ for $i\in [2]$.  This completes the proof for another direction and the base step that $k=2$.

Now we prove the induction step. Let $D$ be a connected acyclic digraph, and let $D'$ be obtained from $D$ by adding two vertices $x,y$ and two arcs $xy, yu$, where $u\in V(D)$. Clearly, $D'$ is also connected and acyclic. We claim that there is a perfect out-forest of size $|V(D)|-k$ in $D$ if and only if there is a perfect out-forest of size $|V(D')|-(k+1)$ in $D'$. Indeed, If $D$ contains a perfect out-forest $F$, then $F$ consists of $k$ disjoint out-trees, that is, $F=\{T_i\mid i\in [k]\}$. So $F'=F\cup \{xy\}$ is a perfect out-forest in $D'$ of size $|V(D')|-(k+1)$. Conversely, If $D'$ contains a perfect out-forest $F'$ of size $|V(D')|-(k+1)$, then $\{xy\}\in F'$ (If $\{xy\}\not\in F'$, then there is an out-tree $T\in F'$ containing $xy$ and $yu$, but this means that the degree of $y$ in $T$ is even, a contradiction). Now observe that $F=F'\setminus \{xy\}$ is a perfect out-forest in $D$ of size $|V(D)|-k$.  

By induction on $k$, we can prove that the problem of deciding whether a given connected ayclic digraph $D$ of order $n$ contains a perfect out-forest with at least $n-k$ arcs is NP-hard, where $k\geq 2$. This completes the proof of the theorem.
\end{pf}

Similarly, the problem of finding a 1-perfect out-forest of size at least $n-1$ is also polynomial-time solvable because $D$ has a 1-perfect forest of size at least $n-1$ if
and only if $D$ is an out-branching in which exactly one vertex is of even degree in $UG(D)$. However, as shown in Theorem~\ref{thmC}, when ``$n-1$'' is replaced by ``$n-k$'' ($k\geq 2$), the problem above becomes NP-hard.

\begin{thm}\label{thmC}
Let $k\geq 2$ be an integer. The problem of deciding whether a given connected digraph $D$ of order $n$ contains a 1-perfect out-forest with at least $n-k$ arcs is NP-hard.
\end{thm}
\begin{pf}
We first consider the case that $k=2$. Let $D'$ be a connected digraph constructed from the digraph $D$ in Theorem~\ref{thmB} such that $V(D')=V(D)\cup \{x\}$ and $A(D')=A(D)\cup \{xx_1^1\}$. With a similar argument to that of Theorem~\ref{thmB}, we can show that $D'$ contains a 1-perfect forest of size at least $n-2$ if and only if $I$ is NAE-satisfied. By the NP-hardness of NAE 3-SAT problem, the problem of deciding whether a given connected digraph of order $n$ contains a 1-perfect out-forest with at least $n-2$ arcs is NP-hard.

We next consider the case that $k\geq 3$. Let $D$ be a connected digraph, we construct another connected digraph $D'$ from $D$ by adding a new vertex $v$ and two arcs $uv, vu$ for some vertex $u\in V(D)$. If $D$ contains a perfect out-forest $F$ of size $|V(D)|-k+1$, then $F$ consists of $k-1$ disjoint out-trees, that is, $F=\{T_i\mid i\in [k-1]\}$. So $F'=F\cup \{v\}$ is a 1-perfect out-forest in $D'$ of size $|V(D')|-k$, where $v$ is the unique vertex whose degree in the underlying graph of $F'$ is even. If $D'$ contains a 1-perfect out-forest $F'$ of size $|V(D')|-k$, then $\{v\}\in F'$. Indeed, if $\{v\}\not\in F'$, then there is an out-tree $T\in F'$ containing $u, v$. As $T$ is an induced subgraph of $D$, we have $uv, vu\in A(T)$, a contradiction. Hence $\{v\}\in F'$ and $v$ is the unique vertex whose degree in the underlying graph of $F'$ is even. Now observe that $F=F'\setminus \{v\}$ is a perfect out-forest in $D$ of size $|V(D)|-k+1$. Therefore, $D$ contains a perfect out-forest of size $|V(D)|-k+1$ if and only if $D'$ contains a 1-perfect out-forest of size $|V(D')|-k$, by Theorem~\ref{thmB}, the result holds.
\end{pf}

Theorems~\ref{thmB} and~\ref{thmC} directly imply the following result.

\begin{thm}\label{thmD}
The following assertions hold:
\begin{description}
\item[(a) ] It is NP-hard to find a perfect out-forest of maximum size in a connected acyclic digraph.
\item[(b) ] It is NP-hard to find an 1-perfect out-forest of maximum size in a connected digraph.
\end{description}
\end{thm}

However, when restricted to some digraph clasess, such as semicomplete digraphs, the above two problems become polynomial-time solvable. As shown in the proof of Theorem~\ref{thmA}(b), $F$ is a perfect out-forest of a semicomplete digraph $D$ if and only if $F$ is a perfect matching of $D$ such that if $xy\in A(F)$, then $yx\not \in A(D)$. 
Hence, each perfect out-forest of $D$ has $|V(D)|/2$ arcs, so, with a similar argument to that of Theorem~\ref{thmA}(b), we can in polynomial-time find a perfect out-forest (if exists) of maximum size in a semicomplete digraph. Similarly, we can also in polynomial-time find a 1-perfect out-forest (if exists) of maximum size in a semicomplete digraph.




\section{Steiner cycle packing}

The problem of {\sc Directed $k$-Linkage} is defined as follows: Given a digraph $D$ and a (terminal) sequence $(s_1, t_1, \dots, s_k, t_k)$ of  distinct vertices of $D,$
decide whether $D$ has $k$ vertex-disjoint paths $P_1, \dots, P_k$, where $P_i$ starts at $s_i$ and ends at $t_i$ for all $i\in [k].$

Sun and Yeo proved the following NP-completeness of {\sc Directed 2-Linkage} for Eulerian digraphs.

\begin{thm}\label{thm101}\cite{Sun-Yeo}
The problem of {\sc Directed 2-Linkage} restricted to Eulerian digraphs is NP-complete.
\end{thm}




Using Theorem~\ref{thm101}, we can prove the following NP-completeness of deciding whether $\kappa^c_S(D)\geq \ell$ for Eulerian digraphs (and therefore for general digraphs).

\begin{thm}\label{thm1a}
Let $k\geq 2, \ell \geq 1$ be fixed integers. For any Eulerian digraph $D$ and $S
\subseteq V(D)$ with $|S|=k$, deciding whether $\kappa^c_S(D) \geq \ell$ is NP-complete.
\end{thm}
\begin{pf}
It is not difficult to see that the problem belongs to NP.
We will show that the problem is NP-hard by reducing from {\sc Directed 2-Linkage} in Eulerian digraphs.
Let $[H; s_1,s_2,t_1,t_2]$ be an instance of {\sc Directed 2-Linkage} in Eulerian digraphs, that is, $H$ is an Eulerian digraph, and $(s_1, t_1, s_2, t_2)$ is a (terminal) sequence of  distinct vertices of $H$. 

We first consider the case that $k\geq 2, \ell\geq 2$. We now produce a new Eulerian digraph $D$ as follows.

Let $$V(H') = V(H) \cup S \cup \{r_1, r_2\},$$
where $S=\{x_i \mid i\in [k]\}$, and let 
\[
\begin{array}{rcl}
 A(H') & = & A(H) \cup \{x_{k-1}s_1, t_1x_k, x_ks_2, t_2x_1, s_1r_1, r_1t_2, s_2r_2, r_2t_1\} \\
        &   & \cup \; \{ x_ix_{i+1} \mid i\in [k] \},\\
\end{array}
\]
where $x_{k+1}=x_1 $. Furthermore, we duplicate the arc $x_ix_{i+1}$ $\ell-1$ times for each $i\in [k-2]$ and duplicate the two arcs $x_{k-1}x_k$ and $x_kx_1$ $\ell-2$ times (note that when $k=2$, we just duplicate the two arcs $x_1x_2$ and $x_2x_1$ $\ell-2$ times). Finally, to avoid parallel arcs, insert a new vertex $z_{i,i+1}^j$ to each arc of the form $x_ix_{i+1}$, where $j\in [\ell]$ if $i\in [k-2]$ and $j\in [\ell-1]$ otherwise. The resulting digraph is called $D$ as shown in Figure~\ref{figure1}.


	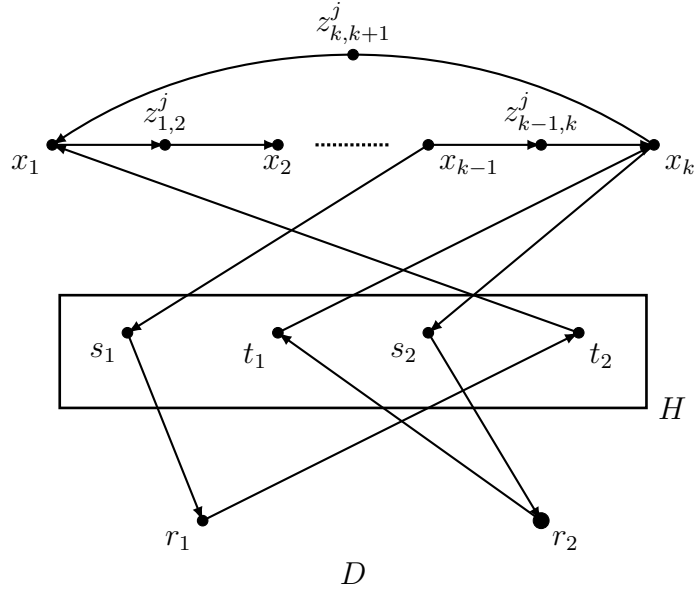
\begin{figure}[htb]
		\centering
		\begin{tikzpicture}
			\filldraw[black]    (0, 0)  circle (2pt)  node [anchor=north east] {$x_{1}$};
			\filldraw[black]    (1.5, 0)  circle (2pt)  node [anchor=south] {$z^{j}_{1,2}$};
			\filldraw[black]    (3, 0)  circle (2pt)  node [anchor=north] {$x_{2}$};
			\filldraw[black]    (5, 0)  circle (2pt)  node [anchor=north west] {$x_{k-1}$};
			\filldraw[black]    (6.5, 0)  circle (2pt)  node [anchor=south] {$z^{j}_{k-1,k}$};
			\filldraw[black]    (8, 0)  circle (2pt)  node [anchor=north west] {$x_{k}$};
			\filldraw[black]    (4, 1.2)  circle (2pt)  node [anchor=south] {$z^{j}_{k,k+1}$};
			\filldraw[black]    (1, -2.5)  circle (2pt)  node [anchor=north east] {$s_{1}$};
			\filldraw[black]    (3, -2.5)  circle (2pt)  node [anchor=north east] {$t_{1}$};
			\filldraw[black]    (5, -2.5)  circle (2pt)  node [anchor=north east] {$s_{2}$};
			\filldraw[black]    (7, -2.5)  circle (2pt)  node [anchor=north west] {$t_{2}$};
			\filldraw[black]    (2, -5)  circle (2pt)  node [anchor=north east] {$r_{1}$};
			\filldraw[black]    (6.5, -5)  circle (3pt)  node [anchor=north west] {$r_{2}$};
			
			\draw[line width=1pt] (0.1, -2) rectangle (7.9, -3.5)  node [anchor=west] {$H$};
			
			\draw[densely dotted] [line width=1.2pt]      (3.5, 0) -- (4.5, 0);
			
			\draw[thick] [-latex]      (0, 0) -- (1.5, 0);
			\draw[thick] [-latex]      (1.5, 0) -- (3, 0);
			\draw[thick] [-latex]      (5, 0) -- (6.5, 0);
			\draw[thick] [-latex]      (6.5, 0) -- (8, 0);
			\draw[thick] [-latex]      (7, -2.5) -- (0, 0);
			\draw[thick] [-latex]      (5, 0) -- (1, -2.5);
			\draw[thick] [-latex]      (3, -2.5) -- (8, 0);
			\draw[thick] [-latex]      (8, 0) -- (5, -2.5);
			\draw[thick] [-latex]      (1, -2.5) -- (2, -5);
			\draw[thick] [-latex]      (2, -5) -- (7, -2.5);
			\draw[thick] [-latex]      (6.5, -5) -- (3, -2.5);
			\draw[thick] [-latex]      (5, -2.5) -- (6.5, -5);
			\draw[thick] []      (8, 0) .. controls (6.5, 1) and (5, 1.2) .. (4, 1.2);
			\draw[thick] [-latex]      (4, 1.2) .. controls (3, 1.2) and (1.5, 1) .. (0, 0);
			
			\node at(4, -5.7) {$D$};
		\end{tikzpicture}
		\caption{The digraph $D$.}
		\label{figure1}
	\end{figure}

It can be checked that $D$ is Eulerian. We will show that $\kappa^c_{S}(D)\ge \ell$ if and only if there exist two vertex-disjoint paths, $P_1$ and $P_2$, in $H$ such that $P_i$ is an $(s_i,t_i)$-path, for $i\in [2]$. This will complete the proof of this case by Theorem~\ref{thm101}.

First assume that there exist two vertex-disjoint paths, $P_1$ and $P_2$ in $H$, such that $P_i$ is an $(s_i,t_i)$-path, for $i\in [2]$.
Add the arcs $x_{k-1}s_1, t_1x_k, x_ks_2, t_2x_1$, and the path $x_1, z^{\ell}_{1,2}, x_2, \dots, z^{\ell}_{k-2, k-1}, x_{k-1}$ to $P_1, P_2$ and call the resulting cycle for $C_{\ell}$. Let $C_j$ be the cycle $x_1, z^{j}_{1,2}, x_2, \dots, z^{j}_{k-1, k}, x_{k}, z^{j}_{k, k+1}, x_1$, for each $j\in [\ell-1]$.  It can be checked that the above cycles $C_1, C_2,\ldots, C_{\ell}$ are internally disjoint, which implies that  $\kappa^c_S(D) \geq \ell$ as desired.

Conversely, assume that $\kappa^c_S(D) \geq \ell$, that is, there is a set of internally disjoint $S$-cycles, say $\{C_i\mid i\in [\ell]\}$. By the construction of $D$, the subdigraph $D[S]$ is the union of the path $x_1, z^{\ell}_{1,2}, x_2, \dots, z^{\ell}_{k-2, k-1}, x_{k-1}$ and $\ell-1$ arc-disjoint cycles: $x_1, z^{j}_{1,2}, x_2, \dots, z^{j}_{k-1, k}, x_{k}, z^{j}_{k, k+1}, x_1$ where $j\in [\ell-1]$. Since $deg_{D}^+(x_i)=\ell$, each $C_i$ contains precisely one out-neighbour and one in-neighbour of $x_i$. There must be two paths, say $P_1, P_2$, such that $P_1$ contains $s_1, t_1$ (resp. $P_2$ contains $s_2, t_2$), and starts from $x_{k-1}$ (resp. $x_k$) and ends at $x_k$ (resp. $x_{k-1}$) in $C_i$ (resp. $C_j$) for some $i\in [\ell]$ (resp. $j\in [\ell]$).

If $i=j$, then clearly there is a pair of vertex-disjoint $s_1-t_1$ path and $s_2-t_2$ path in the cycle $C_i$ (and therefore in $H$), as desired. It remains to consider the case that $i\neq j$. Without loss of generality, assume that $i=1, j=2$. Observe that there is an $x_k-x_{k-1}$ path $P'_1$ in $C_1$ and an $x_{k-1}-x_k$ path $P'_2$ in $C_2$, that is, $C_i=P_i\cup P'_i$ for each $i\in [2]$. Let $C'_1=P_1\cup P_2$ and $C'_2=P'_1\cup P'_2$. It can be checked that $\{C'_i, C_j\mid 1\leq i\leq 2, 3\leq j\leq \ell\}$ is a set of internally disjoint $S$-cycles. With a similar argument to that of the above case, we still have that
there is pair of vertex-disjoint $s_1-t_1$ path and $s_2-t_2$ path in the cycle $C_1$ (and therefore in $H$).

We next consider the case that $k\geq 2, \ell=1$. Similar to the constructon of $H'$ in the case that $k\geq 2, \ell\geq 2$, we now produce a new Eulerian digraph $Q$ as follows. Let $$V(Q) = V(H) \cup S \cup \{r_1, r_2\},$$
where $S=\{x_i \mid i\in [k]\}$, and let 
\[
\begin{array}{rcl}
 A(Q) & = & A(H) \cup \{x_{k-1}s_1, t_1x_k, x_ks_2, t_2x_1, s_1r_1, r_1t_2, s_2r_2, r_2t_1\} \\
        &   & \cup \; \{ x_ix_{i+1} \mid i\in [k-2] \}.\\
\end{array}
\]

It can be checked that $Q$ is Eulerian. With a similar but simpler argument to the above case (and so we omit the details), we can show that $\kappa^c_{S}(Q)\ge 1$ if and only if there exist two vertex-disjoint paths, $P_1$ and $P_2$, in $H$ such that $P_i$ is an $(s_i,t_i)$-path, for $i\in [2]$. This completes the proof of this case by Theorem~\ref{thm101}. 
\end{pf}

We can also prove the NP-completeness of deciding whether $\lambda^c_S(D)\geq \ell$ for general digraphs.

\begin{thm}\label{thm1d}
Let $k\geq 2, \ell \geq 1$ be fixed integers. For a digraph $D$ and $S
\subseteq V(D)$ with $|S|=k$, deciding whether $\lambda^c_S(D) \geq \ell$ is NP-complete.
\end{thm}
\begin{pf}
When $k\geq 2,\ell=1$, by definitions, we clearly have $\lambda^c_S(D) \geq 1$ if and only if 
$\kappa^c_S(D) \geq 1$, so the result holds for this case by Theorem~\ref{thm1a}. Hence, in this following, we assume that $k\geq 2, \ell\ge 2$.

It is not difficult to see that the problem belongs to NP. We will show that the problem is NP-hard by reducing from {\sc Directed 2-Linkage} in Eulerian digraphs.
Let $[H; s_1,s_2,t_1,t_2]$ be an instance of {\sc Directed 2-Linkage} in Eulerian digraphs, that is, $H$ is an Eulerian digraph, and $(s_1, t_1, s_2, t_2)$ is a (terminal) sequence of  distinct vertices of $H$.

We first produce a digraph $D$ the same as that in Theorem~\ref{thm1a}. 
Secondly, we construct a new digraph $D'$ from $D$ as follows: replace every vertex $u$ of $H$ by two vertices $u^-$ and $u^+$ such that $u^-u^+$ is an arc in $D'$ and for every $uv\in A(H)$ add an arc $u^+v^-$ to $D'$. Also, for $z\in S \cup \{r_1, r_2\}$, for every arc $zu$ in $D$ add an arc $zu^-$ to $D'$ and for every arc $uz$ add an arc $u^+z$ to $D'$.

It was proved in Theorem~\ref{thm1a} that $\kappa^c_{S}(D) \geq \ell$ if and only if $[H; s_1,s_2,t_1,t_2]$ is a positive instance of {\sc Directed 2-Linkage} in Eulerian digraphs. Therefore, to prove the theorem, it suffices to show that $\lambda^c_S(D') \geq \ell$ if and only if $\kappa^c_{S}(D) \geq \ell$.

We first assume that $\kappa^c_{S}(D) \geq \ell$, that is, there is a set of internally disjoint $S$-cycles in $D$, say $\{C_i\mid i\in [\ell]\}$. By the construction of $D$, we have $deg_{D}^+(x_i)=\ell$, each $C_i$ contains precisely one out-neighbour and one in-neighbour of $x_i$. There must be two paths, say $P_1, P_2$, such that $P_1$ contains $s_1, t_1$ (resp. $P_2$ contains $s_2, t_2$), and starts from $x_{k-1}$ (resp. $x_k$) and ends at $x_k$ (resp. $x_{k-1}$) in $C_i$ (resp. $C_j$) for some $i\in [\ell]$ (resp. $j\in [\ell]$). Now we obtain a cycle $C'_i$ (resp. $C'_j$) in $D'$ from $C_i$ (resp. $C_j$) as follows: replace every vertex $u$ of $V(H)\cap V(C_i)$ (resp. $V(H)\cap C_j$) by two vertices $u^-$ and $u^+$ such that $u^-u^+$ is an arc in $C'_i$ (resp. $C'_j$) and for every $uv\in A(H)\cap A(C_i)$ (resp. $uv\in A(H)\cap A(C_j)$) add an arc $u^+v^-$ to $A(C'_i)$ (resp. $A(C'_j)$). Also, for each $z\in S \cup \{r_1, r_2\}$, for every arc $zu$ in $A(D)\cap A(C_i)$ (resp. $A(D)\cap A(C_j)$) add an arc $zu^-$ to $C'_i$ (resp. $C'_j$), and for every arc $uz$ add an arc $u^+z$ to $C'_i$ (resp. $C'_j$). Then combining with the remaining $S$-cycles in $\{C_i\mid i\in [\ell]\}$, we obtain a set of arc-disjoint $S$-cycles in $D'$, therefore $\lambda^c_S(D') \geq \ell$.

Conversely, we assume that $\lambda^c_S(D') \geq \ell$, that is, there is a set of arc-disjoint $S$-cycles in $D'$, say $\{C_i\mid i\in [\ell]\}$. By the construction of $D'$, we have $deg_{D'}^+(x_i)=\ell$, each $C_i$ contains precisely one out-neighbour and one in-neighbour of $x_i$. There must be two paths, say $Q_1, Q_2$, such that $Q_1$ contains $s^-_1, t^+_1$ (resp. $Q_2$ contains $s^-_2, t^+_2$), and starts from $x_{k-1}$ (resp. $x_k$) and ends at $x_k$ (resp. $x_{k-1}$) in $C_i$ (resp. $C_j$) for some $i\in [\ell]$ (resp. $j\in [\ell]$). Observe that $Q_1$ must be of the form $x_{k-1}, s^-_1, s^+_1, a^-_1, a^+_1, \dots, a^-_p, a^+_p, t^-_1, t^+_1, x_k$ and $Q_2$ must be of the form $x_{k}, s^-_2, s^+_2, b^-_1, b^+_1, \dots, b^-_q, b^+_q, t^-_2, t^+_2, x_{k-1}$ and furthermore, we have $\{a_i\mid i\in [p]\}\cap \{b_j\mid j\in [q]\}=\emptyset$. Now we obtain $C'_i$ (resp. $C'_j$) from $C_i$ (resp. $C_j$) by replacing $Q_1$ (resp. $Q_2$) with the path $x_{k-1}, s_1, a_1, \dots, a_p, t_1, x_k$ (resp. $x_k, s_2, b_1, \dots, b_p, t_2, x_{k-1}$). It can be checked that, combining with the remaining $S$-cycles in $\{C_i\mid i\in [\ell]\}$, we get a set of $\ell$ internally disjoint $S$-cycles in $D$, therefore $\kappa^c_S(D) \geq \ell$.
\end{pf}

Recall that in Theorem~\ref{thm1a}, we showed that when $D$ is an Eulerian digraph, the problem of deciding whether $\kappa^c_S(D)\geq \ell$ with $|S|=k$ is NP-complete, where both $k\geq 2, \ell\geq 1$ are fixed integers. However, when we consider the class of symmetric digraphs, the problem becomes polynomial-time solvable. We start with the following result by Sun and Yeo.

\begin{lem}\label{thmsym}\cite{Sun-Yeo}
Let $D$ be a symmetric digraph and let $s_1,s_2,\ldots,s_r, t_1,t_2,\ldots,t_r$ be vertices in $D$ (not necessarily disjoint) and
let $S \subseteq V(D)$. We can in $O(|V(G)|^3)$ time decide if there for all $i=1,2,\ldots,r$ exists an $(s_i,t_i)$-path, $P_i$,
such that no internal vertex of any $P_i$ belongs to $S$ or to any path $P_j$ with $j \not=i$ (the end-points of $P_j$ can also not
be internal vertices of $P_i$).
\end{lem}

By Lemma~\ref{thmsym}, 
we will now prove the polynomiality for $\kappa^c_{S}(D)$ on symmetric digraphs.

\begin{thm}\label{thm1c}
Let $k\geq 2$ and $\ell \geq 1$ be fixed integers. We can in polynomial time decide if $\kappa^c_{S}(D) \geq \ell$ for any symmetric
digraph $D$ with $S \subseteq V(D)$, where $|S|=k$.
\end{thm}
\begin{pf}
Let $k\geq 2$ and $\ell \geq 1$ be fixed integers, and let $D$ be a symmetric digraph. Let $S \subseteq V(D)$ with $|S|=k$, and let $\{A_i\mid 0\leq i\leq \ell\}$ be a partition of $A[S]$, where $A[S]$ denotes the arc set in $D[S]$. We need to prove the following claim.

\textbf{Claim:} 
We can in time $O(n^3k^{k\ell})$ decide if there exists a set of $\ell$ internally disjoint $S$-cycles, say $\{C_i \mid i\in [{\ell}]\}$, such that $A(C_i) \cap A[S] = A_i$ for each $i\in [\ell]$ (note that $A_0$ are the arcs in $A[S]$ not used in any of the cycles).

\textbf{Proof of the claim:} Let $C$ be any $S$-cycle in $D$. We define the {\em skeleton} of $C$ as the cycle we obtain from $T$ by contracting all vertices in $V(C)\setminus S$.
Let $C^s$ be a skeleton of an $S$-cycle in $D$. Note that $V(C^s)=S$ and therefore there  are at most $k^k$ different skeletons of $S$-cycles in $D$.

Our algorithm will try all possible $\ell$-tuples, ${\cal C}^s = (C_1^s, C_2^s, \ldots, C_{\ell}^s)$, of skeletons of $S$-cycles and determine if there is a set of $\ell$ internally disjoint $S$-cycles, say $\{C_i \mid i\in [{\ell}]\}$, such that $A(C_i) \cap A[S] = A_i$ and  $C_i^s$ is the skeleton of $C_i$ for each $i\in [\ell]$. If such a set of cycles exists for any ${\cal C}^s$, then we return this solution, and if no such set of cycles exist for any ${\cal C}^s$, then we return that no solution exists. We will prove that this algorithm gives the correct answer and compute its time complexity.

If our algorithm returns a solution, then clearly a solution exists. So now assume that a solution exists and let $\{C_i \mid i\in [{\ell}]\}$ be the desired set of internally disjoint $S$-cycles. When we consider ${\cal C}^s = (C_1^s, C_2^s, \ldots, C_{\ell}^s)$, where $C_i^s$ is the skeleton of $C_i$, our algorithm will find a
solution, so the algorithm always returns a solution if one exists.


Given such an $\ell$-tuples, ${\cal C}^s$, we need to determine if there is a set of $\ell$ internally disjoint $S$-cycles, say $\{C_i \mid i\in [{\ell}]\}$, such that $A(C_i) \cap A[S] = A_i$ and  $C_i^s$ is the skeleton of $C_i$ for each $i\in [\ell]$.
We first check that the arcs in $A_i$ belong to the skeleton $C_i^s$ and that no vertex in $V(D)\setminus S$ belongs to more than one skeleton.
If the above does not hold, then the desired cycles do not exist, so we assume that the above holds in the following argument. For every arc $uv \not\in A[S]$ that belongs to some skeleton $C_i^s$, we want to find a $(u,v)$-path in $D - A[S]$,
such that no internal vertex on any path belongs to $S$ or to a different path.
This can be done in $O(n^3)$ time by Lemma~\ref{thmsym}.  If such paths exist, then we obtain the desired $S$-cycles by substituting each $uv$ by the corresponding $(u,v)$-path. Otherwise, the desired set of $S$-cycles does not exist. 

Therefore, the algorithm works correctly and has complexity $O(n^3k^{k\ell})$ since the number of different $\ell$-tuples, ${\cal C}^s$, that we need to consider is bounded by the function $(k^k)^{\ell}=k^{k\ell}$.\qed

By the claim, 
we can in time $O(n^3k^{k\ell})$ decide if there exists a set of $\ell$ internally disjoint $S$-cycles, say $\{C_i \mid i\in [{\ell}]\}$, such that $A(C_i) \cap A[S] = A_i$ for each $i\in [\ell]$.
We will now use the algorithm of the claim for all possible partitions ${\cal A} = \{A_i\mid 0\leq i\leq \ell\}$. If we find the desired set of $\ell$ internally disjoint $S$-cycles for any such a partition, then we return ``$\kappa^c_{S}(D) \geq \ell$''. Otherwise, we return ``$\kappa^c_{S}(D) < \ell$''.
Note that if  $\kappa^c_{S}(D) \geq \ell$, then we will correctly determine that $\kappa^c_{S}(D) \geq \ell$, when we consider the correct partition ${\cal A}$, which proves that the above algorithms will always return the correct answer.

Hence, we deduce that the algorithm works correctly and has complexity $O(n^3k^{k\ell}(\ell+1)^{k^2/2})$ which is a polynomial in $n$, since the number of partitions ${\cal A}$ of $A(D[S])$ is bounded by $(\ell+1)^{|A(D[S])|} \leq (\ell+1)^{k^2/2}$, and the two parameters $k$ and $\ell$ are fixed.
\end{pf}

By Theorem~\ref{thm1c}, we directly have the following result.

\begin{cor}\label{cor11}
We can in polynomial time decide if a given symmetric digraph $D$ is $k$-cyclic, for a fixed integer $k\geq 2$.
\end{cor}

\section{Discussions}

The perfect forest theorem in undirected graphs shows that every connected graph of even order contains a 0-perfect forest \cite{Caro-Lauri-Zarb, GutinJGT2016, ScottGC2001}. However, there is no similar result on $i$-perfect out-forest for digraphs ($i\in \{0,1\}$), that is, not every strong connected digraph of even order has a perfect out-forest, and not every strong connected digraph of odd order has a 1-perfect out-forest, as any symmetric digraph $D$ contains neither perfect out-forest nor 1-perfect out-forest. Hence, it would be interesting to give some nice sufficient conditions or characterizations to guarantee the existence of an $i$-perfect out-forest in digraphs ($i\in \{0,1\}$).

Note that the construction in Theorem~\ref{thm1d} may not make $D'$ Eulerian even when $H$ is Eulerian, since $deg^-_{D'}(u^-)(=deg^-_{D}(u))$ maybe larger than $deg^+_{D'}(u^-)(=1)$. Therefore, we need to find other approach to prove the NP-completeness for $\lambda^c_S(D)$ for an Eulerian digraph $D$. Certainly, we cannot exclude the possibility that deciding whether $\lambda^c_S(D)\geq \ell$ restricted to Eulerian digraphs is polynomial-time solvable.






\vskip 1cm

\noindent {\bf Data Availability Statement.} Data sharing not applicable to this article
as no datasets were generated or analysed during the current study.

\vskip 1cm

\noindent {\bf Acknowledgement.} We are thankful to Professor Anders Yeo for discussions on the perfect forests problem. Yuefang Sun was supported by Yongjiang Talent Introduction
Programme of Ningbo under Grant Number 2021B-011-G and Zhejiang Provincial Natural Science Foundation of China under Grant Number LY20A010013. 

\end{document}